\documentclass{amsart}

\newtheorem{theorem}{Theorem}[section]
\newtheorem{lemma}[theorem]{Lemma}

\newtheorem{proposition}[theorem]{Proposition}
\newtheorem{corollary}[theorem]{Corollary}

\theoremstyle{definition}
\newtheorem{definition}[theorem]{Definition}
\newtheorem{example}[theorem]{Example}

\theoremstyle{remark}
\newtheorem{remark}[theorem]{Remark}

\usepackage{graphicx}
\usepackage{color}
\usepackage{amsmath}
\usepackage{amsfonts}
\usepackage{amssymb}
\newcommand{\be}{\begin{equation}}
\newcommand{\ee}{\end{equation}}

\newcommand{\al}{\alpha}

\newcommand{\bet}{\beta}

\newcommand{\g}{\gamma}

\newcommand{\Om}{\Omega}

\newcommand{\om}{\omega}

\newcommand{\we}{{\stackrel{\scriptscriptstyle{W}}{\Gamma}}\phantom{}}
\newcommand{\wc}{{\stackrel{\scriptscriptstyle{C}}{\Gamma}}\phantom{}}
\newcommand{\rw}{{\stackrel{\scriptscriptstyle{W}}{R}}\phantom{}}
\newcommand{\omw}{{\stackrel{\scriptscriptstyle{W}}{\Omega}}\phantom{}}

\newcommand{\con}{{\stackrel{\scriptscriptstyle{0}}{\Gamma}}\phantom{}}

\newcommand{\dz}{\wedge}

\newcommand{\ba}{\begin{array}}

\newcommand{\ea}{\end{array}}

\newcommand{\beq}{\begin{eqnarray}}

\newcommand{\eeq}{\end{eqnarray}}

\newtheorem{lm}{lemma}

\newtheorem{thee}{theorem}

\newtheorem{proo}{proposition}

\newtheorem{co}{corollary}

\newtheorem{rem}{remark}

\newtheorem{deff}{definition}

\newcommand{\bd}{\begin{deff}}

\newcommand{\ed}{\end{deff}}

\newcommand{\bl}{\begin{lm}}

\newcommand{\el}{\end{lm}}

\newcommand{\bp}{\begin{proo}}

\newcommand{\ep}{\end{proo}}

\newcommand{\bt}{\begin{thee}}

\newcommand{\et}{\end{thee}}

\newcommand{\bc}{\begin{co}}

\newcommand{\ec}{\end{co}}

\newcommand{\brm}{\begin{rem}}

\newcommand{\erm}{\end{rem}}

\newcommand{\der}{{\rm d}}

\newcommand{\sgn}{\mathrm{sgn}}

\usepackage{t1enc}
\def\frak{\mathfrak}

\newcommand{\newc}{\newcommand}

\newcommand{\x}{\times}

\let\ccdot\cdot
\def\cdot{\hbox to 2.5pt{\hss$\ccdot$\hss}}

\newc{\aR}{\mbox{\boldmath{$ R$}}}
\newc{\aS}{\mbox{\boldmath{$ S$}}}
\newc{\aT}{\mbox{\boldmath{$ T$}}}
\newc{\aW}{\mbox{\boldmath{$ W$}}}

\newc{\aK}{\mbox{\boldmath{$ K$}}}
\newc{\aL}{\mbox{\boldmath{$ L$}}}

\usepackage{amssymb}
\usepackage{amscd}

\let\i=\iota

\newcommand{\hook}{\raisebox{-0.35ex}{\makebox[0.6em][r]
{\scriptsize $-$}}\hspace{-0.15em}\raisebox{0.25ex}{\makebox[0.4em][l]{\tiny
 $|$}}}

\newc{\obstrn}[2]{B^{#1}_{#2}}

\newcommand{\rpl}                         
{\mbox{$
\begin{picture}(12.7,8)(-.5,-1)
\put(0,0.2){$+$}
\put(4.2,2.8){\oval(8,8)[r]}
\end{picture}$}}

\newcommand{\lpl}                         
{\mbox{$
\begin{picture}(12.7,8)(-.5,-1)
\put(2,0.2){$+$}
\put(6.2,2.8){\oval(8,8)[l]}
\end{picture}$}}

\usepackage{ifthen}

\newcommand{\bbR}{\mathbb{R}}
\newcommand{\bbN}{\mathbb{N}}
\newcommand{\bfa}{\mathbf{a}}
\newcommand{\sog}{\mathbf{SO}}
\newcommand{\cog}{\mathbf{CO}}

\newcommand{\slg}{\mathbf{SL}}
\newcommand{\glg}{\mathbf{GL}}
\newcommand{\og}{\mathbf{O}}
\newcommand{\coa}{\frak{co}}
\newcommand{\soa}{\frak{so}}
\newcommand{\sla}{\frak{sl}}
\newcommand{\gla}{\frak{gl}}

\newcommand{\spg}{\mathbf{Sp}}
\newcommand{\sug}{\mathbf{SU}}

\newc{\tensor}[1]{#1}
\newc{\Mvariable}[1]{\mbox{#1}}
\newc{\down}[1]{{}_{#1}}
\newc{\up}[1]{{}^{#1}}

\newc{\JulyStrut}{\rule{0mm}{6mm}}
\newc{\midtenPan}{\mbox{\sf S}}
\newc{\midten}{\mbox{\sf T}}
\newc{\midtenEi}{\mbox{\sf U}}
\newc{\ATen}{\mbox{\sf E}}
\newc{\BTen}{\mbox{\sf F}}
\newc{\CTen}{\mbox{\sf G}}

\def\sideremark#1{\ifvmode\leavevmode\fi\vadjust{\vbox to0pt{\vss
 \hbox to 0pt{\hskip\hsize\hskip1em
 \vbox{\hsize3cm\tiny\raggedright\pretolerance10000
 \noindent #1\hfill}\hss}\vbox to8pt{\vfil}\vss}}}%

\newcommand{\edz}[1]{\sideremark{#1}}
\newcommand{\bgw}{{\textstyle \bigwedge}}
\newcommand{\bgs}{{\textstyle \bigodot}}
\newcommand{\bgt}{{\textstyle \bigotimes}}

\numberwithin{equation}{section}

\newcommand{\bma}{\begin{pmatrix}}
\newcommand{\ema}{\end{pmatrix}}

\newcounter{romenumi}
\newcommand{\labelromenumi}{(\roman{romenumi})}

\newcommand{\ten}{\Upsilon}

\begin{document}
\title{Comment on ${\bf GL}(2,\bbR)$ geometry of 4th order ODE's }

\author{Pawe\l~ Nurowski} 
\address{Instytut Fizyki Teoretycznej,
Uniwersytet Warszawski, ul. Hoza 69, Warszawa, Poland}
\email{nurowski@fuw.edu.pl} 
\thanks{This research was supported by the KBN grant 1 P03B 07529}
\date{\today}

\begin{abstract}
We describe 4th order ODEs satisfying two contact invariant conditions
of Bryant in terms of the Ricci tensor of a certain 
$\gla(2,\bbR)$ valued connection. We also provide nonhomogeneous
examples of such ODEs.
\vskip5pt\centerline{\small\textbf{MSC classification}: 53A40, 53B05,
  34C30, 34C31}\vskip15pt
\end{abstract}
\maketitle
\tableofcontents
\section{Introduction}
Recently there is a growing interest in the geometrization program of
ODEs \cite{dub,dun,gn0,gn,nursp,nur}. Although the program my be traced back to
S. Lie \cite{Lie} and M. A. Tresse \cite{tres}, and although it was
formulated by E. Cartan and S. S. Chern in the 1940s
\cite{Cartpc,CartanSpanish,chern}, it was not very popular until the
works of R. Bryant (see e.g. \cite{bryantsp4}) on the invariants of
the fourth order ODEs. In the present note we restate some of the
results of \cite{bryantsp4} in terms of the invariants of the recently discussed
$\glg(2,\bbR)$ geometry of ODEs \cite{gn}. In particular we interpret
Bryant's results in terms of the Ricci tensor of a certain 
$\gla(2,\bbR)$ connection, which characterises the ODEs satisfying
contact invariant conditions of Bryant \cite{bryantsp4}.

Our starting point is the following well known
\begin{proposition}\label{odef2}
Ordinary differential equation
$$y^{(4)}=0$$
has $\glg(2,\bbR)\times_{\rho}\bbR^4$ as its group of contact
symmetries. Here $\rho :\glg(2,\bbR)\to\glg(4,\bbR)$ is the
$4$-dimensional irreducible representation of $\glg(2,\bbR)$. 
\end{proposition}

The representation $\rho$, at the level of the Lie algebra
$\gla(2,\bbR)$, is given in terms of the Lie algebra generators
\begin{eqnarray}
&&E_+=\bma
0&3&0&0\\0&0&2&0\\0&0&0&1\\0&0&0&0\ema,\quad
E_-=\bma
0&0&0&0\\1&0&0&0\\0&2&0&0\\0&0&3&0\ema,\quad
E_0=\bma
-3&0&0&0\\0&-1&0&0\\0&0&1&0\\0&0&0&3\ema,\nonumber\\&&\quad\quad\quad\quad\quad\quad\quad\quad\quad\quad\quad E=-\bma
3&0&0&0\\0&3&0&0\\0&0&3&0\\0&0&0&3\ema.\label{rep}
\end{eqnarray}
These matrices satisfy the $\gla(2,\bbR)$ commutation
relations
$$[E_0,E_+]=-2E_+\quad,\quad [E_0,E_-]=2E_-\quad,\quad
[E_+,E_-]=-E_0\quad,$$
where the commutator in the $\gla(2,\bbR)={\rm
  Span}_\bbR(E_-,E_+,E_0,E)\subset End(\bbR^4)$ is 
the usual commutator of matrices. 
   
Now, we consider a general $4$-th order
ODE
\be
y^{(4)}=F(x,y,y',y'',y^{(3)}).
\label{ode}
\ee
To simplify the notation, we introduce the coordinates 
$x,y,y_1=y',y_2=y'',y_3=y^{(3)}$ on the 
$5$-dimensional \emph{jet space} $J$. Introducing the four \emph{contact forms} 
\begin{eqnarray}
&&\om^0=\der y-y_1\der x\nonumber\\
&&\om^1=\der y_1-y_2\der x\label{fomry}\\
&&\om^{2}=\der y_{2}-y_{3}\der x\nonumber\\
&&\om^{3}=\der y_{3}-F(x,y,y_1,y_2,y_3)\der x,\nonumber
\end{eqnarray}
and an additional 1-form 
$$w_+=\der x,$$
we define a \emph{contact transformation} to be a diffeomorphism $\phi:J\to
J$ which transforms the above five one-forms via:
\begin{eqnarray}
&&\phi^*\om^0=\al^0_{~0}\om^0\nonumber\\
&&\phi^*\om^1=\al^1_{~0}\om^0+\al^1_{~1}\om^1\nonumber\\
&&\phi^*\om^2=\al^2_{~0}\om^0+\al^2_{~1}\om^1+\al^2_{~2}\om^2\label{cont}\\
&&\phi^*\om^3=\al^3_{~0}\om^0+\al^3_{~1}\om^1+\al^3_{~2}\om^2+\al^3_{~3}\om^3\nonumber\\
&&\phi^*w_+=\al^4_{~0}\om^0+\al^4_{~1}\om^1+\al^4_{~4}w_+.\nonumber
\end{eqnarray} 
Here $\al^i_{~j}$, $i,j=0,1,2,3,4,5$, are real functions on $J$ such that
$$\al^0_{~0}\al^1_{~1}\al^2_{~2}\al^3_{~3}\al^4_{~4}\neq 0.$$
The contact equivalence problem for the $4$th order ODEs (\ref{ode}) 
can be studied in terms of the invariant forms 
$(\theta^0,\theta^1,\theta^2\theta^3,\Om_+)$ defined by
\be
\bma\theta^0\\\theta^1\\\theta^2\\\theta^3\\\Om_+\ema
=\bma\al^0_{~0}&&&&&\\
\al^1_{~0}&\al^1_{~1}&&&&\\
\al^2_{~0}&\al^2_{~1}&\al^2_{~2}&&&\\
\al^3_{~0}&\al^3_{~1}&\al^3_{~2}&\al^3_{~3}&&\\
\al^4_{~0}&\al^4_{~1}&&&&\al^4_{~4}\ema
\bma\om^0\\\om^1\\\om^2\\\om^3\\w_+\ema.\label{fomr}
\ee
Among all ODEs (\ref{ode}) considered modulo contact transformations
(\ref{cont}) there is a remarkable class for which the invariant forms
satisfy
\begin{eqnarray}
&&\der\theta^0=3(\Om+\Om_0)\dz\theta^0-3\Om_+\dz\theta^1\nonumber\\
&&\der\theta^1=-\Om_-\dz\theta^0+(3\Om+\Om_0)\dz\theta^1-2\Om_+\dz\theta^2\nonumber\\
&&\der\theta^2=-2\Om_-\dz\theta^1+(3\Om-\Om_0)\dz\theta^2-\Om_+\dz\theta^3\label{ffs}\\
&&\der\theta^3=-3\Om_-\dz\theta^2+3(\Om-\Om_0)\dz\theta^3.\nonumber
\end{eqnarray}
This system is defined on an 8-dimensional $\glg(2,\bbR)$ principal
fibre bundle $P$ over the solution space $M^4$ for the corresponding ODE
(\ref{ode}). The invariant forms
$(\theta^0,\theta^1,\theta^2,\theta^3,\Om_+)$ together with the additional
three 1-forms $(\Om_-,\Om_0,\Om)$ constitute a well defined coframe on
$P$.  

As noted by Bryant \cite{bryantsp4}, the class of ODEs having forms
$(\theta^0,\theta^1,\theta^2,\theta^3,\Om_+,\Om_-,\Om_0,\Om)$ of
system (\ref{ffs}), is distinguished by the demand that their defining
functions $F=F(x,y,y_1,y_2,y_3)$  
satisfy the following two conditions: 
\begin{eqnarray}
&&4D^2F_3-8DF_2+8F_1-6DF_3 F_3+4F_2 F_3+F_3^3=0,\nonumber\\
&&\label{br}\\
&&160D^2F_2-640DF_1+144(DF_3)^2-352 DF_3 F_2+144 F_2^2 -\nonumber\\
&&80DF_2 F_3+160 F_1 F_3-72 DF_3 F_3^2+88 F_2 F_3^2+9 F_3^4+16000F_y=0.\nonumber\end{eqnarray}
Here $F_i=\tfrac{\partial F}{\partial y_i}$ and 
$D=\partial_x+y_1\partial_y+y_2\partial_{y_1}+y_3\partial_{y_2}+F\partial_{y_3}$.
Bryant's conditions (\ref{br}), considered simultaneously, 
are \emph{contact invariant}; if the ODE undergoes contact
transformation of its variables, the conditions (\ref{br}) are
preserved. 
Examples are known of ODEs satisfying these conditions \cite{bryantsp4}, the simplest
being
\be 
y^{(4)}=(y^{(3)})^{(4/3)}.\label{examm}
\ee
The purpose of this note is to establish
a theorem on speciality of a $\gla(2,\bbR)$-valued connection defined
by such ODEs on their solution spaces.
\section{The closed system}
Let us make the following choice 
\begin{eqnarray}
\al^0_{~0}&=&-3\al^1_{~1}\al^4_{~4}\nonumber\\
\al^2_{0}&=&-\frac{(\al^1_{~0})^2}{3\al^1_{~1}\al^4_{~4}}+\frac{\al^1_{~1}}{240\al^4_{~4}}(-24DF_3
  + 36 F_2 + 11 F_3^2)\nonumber\\
\al^2_{~1}&=&-\frac{2\al^1_{~0}}{3\al^4_{~4}}
+\frac{\al^1_{~1}}{12\al^4_{~4}} F_3\nonumber\\
\al^2_{~2}&=&-\frac{\al^1_{~1}}{2\al^4_{~4}}\nonumber\\
\al^3_{~0}&=&\frac{(\al^1_{~0})^3}{9(\al^1_{~1}\al^4_{~4})^2}
+\frac{\al^1_{~0}}{240(\al^4_{~4})^2}(24 DF_3 - 36 F_2 - 11F_3^2) +\nonumber\\&&
\frac{\al^1_{~1}}{720(\al^4_{~4})^2}(36(DF_2 - 4 F_1) 
+ 18(DF_3 - 2 F_2) F_3 - 7 F_3^3)\label{ala}\\
\al^3_{~1}&=&\frac{(\al^1_{~0})^2}{3\al^1_{~1}(\al^4_{~4})^2}
-\frac{\al^1_{~0}}{12(\al^4_{~4})^2}F_3 +
\frac{\al^1_{~1}}{240(\al^4_{~4})^2}(36DF_3 - 84 F_2-19 F_3^2)\nonumber\\
\al^3_{~2}&=&\frac{\al^1_{~0}}{2(\al^4_{~4})^2}-\frac{\al^1_{~1}}{4(\al^4_{~4})^2}F_3\nonumber\\
\al^3_{~3}&=&\frac{\al^1_{~1}}{2(\al^4_{~4})^2}\nonumber\\
\al^4_{~0}&=&-\frac{\al^4_{~4}}{60}(12 DF_{33} - 6 F_{23} + F_3
F_{33})\nonumber\\
\al^4_{~1}&=&\frac{\al^4_{~4}}{6}F_{33}\nonumber
\end{eqnarray}
for the group parameters defining
forms $(\theta^0,\theta^1,\theta^2,\theta^3,\Om_+)$ of
(\ref{fomr}). Then we have the following  
\begin{theorem}\label{prp}
If a fourth order ODE \be
y^{(4)}=F(x,y,y',y'',y^{(3)})\label{ode4}
\ee 
satisfies contact invariant conditions 
\begin{eqnarray}
&&4D^2F_3-8DF_2+8F_1-6DF_3 F_3+4F_2 F_3+F_3^3=0,\nonumber\\
&&\label{wun}\\
&&160D^2F_2-640DF_1+144(DF_3)^2-352 DF_3 F_2+144 F_2^2 -\nonumber\\
&&80DF_2 F_3+160 F_1 F_3-72 DF_3 F_3^2+88 F_2 F_3^2+9
  F_3^4+16000F_y=0\nonumber
\end{eqnarray}
then the manifold $P$ parametrised by
$(x,y,y_1,y_2,y_3,\al^1_{~0},\al^1_{~1},\al^4_{~4})$ is a 
principal $\glg(2,\bbR)$ bundle $P\to M^4$ over the solution
space $M^4$ of (\ref{ode4}) and forms 
$(\theta^0,\theta^1,\theta^2,\theta^3,\Om_+)$, together with
additional three 1-forms $(\Om_-,\Om_0,\Om)$, constitute an invariant 
coframe on $P$ satisfying
\begin{eqnarray}
\der\theta^0&=&3(\Om+\Om_0)\dz\theta^0-3\Om_+\dz\theta^1\nonumber\\
\der\theta^1&=&-\Om_-\dz\theta^0+(3\Om+\Om_0)\dz\theta^1-2\Om_+\dz\theta^2\nonumber\\
\der\theta^2&=&-2\Om_-\dz\theta^1+(3\Om-\Om_0)\dz\theta^2-\Om_+\dz\theta^3\label{sy1}\\
\der\theta^3&=&-3\Om_-\dz\theta^2+3(\Om-\Om_0)\dz\theta^3\nonumber\\&&\nonumber\\
\der\Om_+&=&2\Om_0\dz\Om_++
\tfrac{1}{12}(-3a_0 + 4b_1)\theta^0\dz\theta^1 +\nonumber\\&& 
    \tfrac14(a_1 + 2 b_2)\theta^0\dz\theta^2 + 
    \tfrac{1}{24}(3 a_2 + 4 b_3)\theta^0\dz\theta^3 + \nonumber\\&&
    \tfrac18(-5 a_2 + 4 b_3)\theta^1\dz\theta^2 + 
    \tfrac16 b_4\theta^1\dz\theta^3\nonumber\\
\der\Om_-&=&-2\Om_0\dz\Om_-+
\tfrac16b_0\theta^0\dz\theta^2 +\nonumber\\&& 
    \tfrac{1}{24}(-3a_0 + 4 b_1)\theta^0\dz\theta^3 + 
    \tfrac18 (5a_0 + 4 b_1)\theta^1\dz\theta^2+ \nonumber\\&&
    \tfrac14(-a_1 + 2b_2)\theta^1\dz\theta^3 + 
    \tfrac{1}{12}(3 a_2 + 4 b_3)\theta^2\dz\theta^3\label{sy2}\\
\der\Om_0&=&\Om_+\dz\Om_--
\tfrac16b_0\theta^0\dz\theta^1 +\nonumber\\&& 
    \tfrac{1}{24}(-3 a_0 - 4b_1)\theta^0\dz\theta^2 + 
    \tfrac14a_1\theta^0\dz\theta^3 - \tfrac14 a_1\theta^1\dz\theta^2 +\nonumber\\&& 
    \tfrac{1}{24}(-3 a_2 + 4 b_3)\theta^1\dz\theta^3 + 
    \tfrac16b_4\theta^2\dz\theta^3\nonumber\\
\der\Om&=&-\tfrac16 b_0\theta^0\dz\theta^1 - \tfrac13b_1\theta^0\dz\theta^2 - 
    \tfrac16b_2\theta^0\dz\theta^3 -\nonumber\\&&\tfrac12b_2\theta^1\dz\theta^2 - 
    \tfrac13 b_3\theta^1\dz\theta^3 - \tfrac16b_4\theta^2\dz\theta^3.\nonumber
\end{eqnarray}
The coefficients $a_0$, $a_1$, $a_2$, $b_0$, $b_1$, $b_2$, $b_3$,
$b_4$ are totally determined by (\ref{ode4}) and are expressible in terms of
the derivatives of function $F$ and the coordinates. The simplest of
these coefficients are:
\begin{eqnarray*}
b_4&=&-2\frac{(\al^4_{~4})^3}{(\al^1_{~1})^2}F_{333}
\end{eqnarray*}
\begin{eqnarray*}
b_3&=&-\frac{(\al^4_{~4})^2}{12(\al^1_{~1})^2}~( 6DF_{333}+ 
              5F_3F_{333})
              -\frac{2\al^1_{~0}(\al^4_{~4})^2}{3(\al^1_{~1})^3}F_{333}
\end{eqnarray*}
\begin{eqnarray*}
b_2&=&-\frac{2\al^4_{~4}(\al^1_{~0})^2}{9(\al^1_{~1})^4}F_{333}
  -
\frac{\al^4_{~4}\al^1_{~0}}{18(\al^1_{~1})^3}~(6DF_{333} + 5F_3
F_{333})+\\&&
\frac{\al^4_{~4}}{360(\al^1_{~1})^2}[60(2DF_{233} + 4F_{133} - 2F_{223} + DF_{333}F_3) +(-36 DF_3 + 204F_2 + 79\ F_3^2)F_{333}]
\end{eqnarray*}
\begin{eqnarray*}
a_2&=&-\frac{(\al^4_{~4})^2}{45(\al^1_{~1})^2}~(18 DF_{333} + 24F_{233} + 4 F_{33}^2 + 
              27F_3F_{333}).
\end{eqnarray*}
Other coefficients are given in the next two sections.
\end{theorem}
The proof of this theorem is a lengthy calculation based on a variant
of Cartan's equivalence method. In the next
section we outline the main points of the proof.
\section{Proof of the main theorem}

The basic idea in the proof of Theorem \ref{prp} is to force 1-forms (\ref{fomr}) to
satisfy system (\ref{ffs}). This requirement makes restrictions on the
free parameters $\al^i_{~j}$ and, more importantly, on the possible
functions $F=F(x,y,y',y'',y^{(3)})$ defining the ODE.

The main steps when imposing (\ref{ffs}) on (\ref{fomr}) are:
\begin{itemize} 
\item[1)] equation
  $\der\theta^0\dz\theta^0\dz\theta^2=3\Om_+\dz\theta^0\dz\theta^1\dz\theta^2$
  requires $\al^0_{~0}=-3\al^1_{~1}\al^4_{~4}$,
\item[2)] the first equation (\ref{ffs}) gives a relation between $\Om$,
  $\Om_0$, $\der\al^4_{~4}$ and $\der\al^1_{~1}$,
\item[3)] similarly, equation
  $\der\theta^1\dz\theta^1\dz\theta^2=-\Om_-\dz\theta^0\dz\theta^1\dz\theta^2$
  gives a relation between $\Om_-$, $\der\al^1_{~0}$,
  $\der\al^1_{~1}$,
\item[4)] equation
  $\der\theta^1\dz\theta^0\dz\theta^1=-2\Om_+\dz\theta^0\dz\theta^1\dz\theta^2$
  gives $\al^2_{~2}=-\frac{\al^1_{~1}}{2\al^4_{~4}}$,  
\item[5)] equation
  $\der\theta^1\dz\theta^0\dz\theta^2=-(3\Om+\Om_0)\dz\theta^0\dz\theta^1\dz\theta^2$
  gives a relation between $\Om$, $\Om_0$ and $\der\al^1_{~1}$,
\item[6)] now, the expressions for
  $\der\theta^2\dz\theta^0\dz\theta^1\dz\theta^3$,
  $\der\theta^2\dz\theta^0\dz\theta^1\dz\theta^3$,
  $\der\theta^3\dz\theta^0\dz\theta^1\dz\theta^2$ enable us to fix $\al^3_{~2}$,
  $\al^3_{~3}$, and $\al^2_{~1}$ respectively,
\item[7)] considering successively $\der\theta^3\dz\theta^0\dz\theta^1$,
  $\der\theta^2\dz\theta^0$,
  $\der\theta^2\dz\theta^1\dz\theta^2\dz\theta^3$ we fix $\al^2_{~0}$,
  $\al^3_{~1}$, $\al^3_{~0}$,
\item[8)] now the requirement
  $\der\theta^3\dz\theta^0\dz\theta^2\dz\theta^3=0$ gives the first
  Bryant condition $4D^2F_3-8DF_2+8F_1-6DF_3F_3+4F_2F_3-F_3^3=0$,
\item[9)] the second of Bryant's conditions (\ref{br}) is equivalent to
  the requirement that
  $\der\theta^3\dz\theta^1\dz\theta^2\dz\theta^3=0$,
\item[10)] now, having Braynt's conditions determined, it is
  straightforward to obtain the required system (\ref{ffs}) and 
express all the 
  $\al^i_{~j}$s in terms of $\al^1_{~0}$, $\al^1_{~1}$ and
  $\al^4_{~4}$only,
\item[11)] the expressions for $\al^i_{~j}$s are given by (\ref{ala});
  inserting them to (\ref{fomr}) we get the invariant forms
  $(\theta^0,\theta^1,\theta^2,\theta^3,\Om_+)$
\item[12)] forms $\Om_0,\Om_-,\Om$ are determined by the linear relations
  from points 2), 3) and 5).
\end{itemize}
In this way one finds the explicit expressions for the invariant
coframe satisfying system (\ref{ffs}). Instead of giving these
formulae we present formulae for
$(\theta^0,\theta^1,\theta^2,\theta^3,\Om_+)$ evaluated at 
$(\al^1_{~0},\al^1_{~1},\al^4_{4})=(0,1,1)$. Denoting these forms by 
$(\theta^0_0,\theta^1_0,\theta^2_0,\theta^3_0,\Om^0_+)$, we have

\begin{eqnarray*}
\theta^0_0&=&-3\om^0\\
\theta^1_0&=&\om^1\\
\theta^2_0&=&\tfrac{1}{240}(-24DF_3+36F_2+11F_3^2)\om^0+\tfrac{1}{12}F_3\om^1-\tfrac12\om^2\\
\theta^3_0&=&\tfrac{1}{720}\Big(36(DF_2-4F_1)+18(DF_3-2F_2)F_3-7F_3^2\Big)\om^0+\\
&&\tfrac{1}{240}(36DF_3-84F_2-19F_3^2)\om^1-\tfrac14
F_3\om^2+\tfrac12\om^3\\
\Om^0_+&=&-\tfrac{1}{60}(12DF_{33}-6F_{23}+F_3F_{33})\om^0+\tfrac16 F_{33}\om^1+w_+.
\end{eqnarray*}

The remaining three 1-forms $(\Om_0,\Om_-,\Om)$, when written in the gauge
$(\al^1_{~0},\al^1_{~1},\al^4_{4})=(0,1,1)$, read:
\begin{eqnarray*}
\Om^0_0&=&\frac{1}{4320}(72 DF_{23} + 432F_{13} - 288 F_{22} + 60 DF_{33} F_3 - 216 F_{23} F_3 -\\&& 108DF_3 F_{33} + 324 F_2 F_{33} + 47 F_3^2 F_{33}) \theta^0_0 + 
    \frac{1}{180}(3 DF_{33} - 9 F_{23} - F_3 F_{33})\theta^1_0 +\\&& \frac16 F_{33}\theta^2_0 - 
    \frac{1}{12} F_3\theta^4_0\\
\Om^0_-&=&\frac{1}{64800} (720DF_{22} + 288 DF_3 DF_{33} - 2160 F_{12} - 
              432 DF_{33} F_2 + 216 DF_3 F_{23} +\\&& 216 F_2 F_{23} + 720 DF_{23}F_3 - 1080 F_{13}F_3 - 360 F_{22} F_3 + 48 DF_{33} F_3^2 - \\&&
              174 F_{23}F_3^2 - 360 DF_2 F_{33} + 1440 F_1 F_{33} + 
              24 DF_3 F_3 F_{33} + 324 F_2 F_3 F_{33} +\\&& 29 F_3^3 F_{33} + 
              3600 F_{3y})\theta^0_0 +\\&&\frac{1}{1080}(-108 DF_{23} - 288 F_{13} + 252 F_{22} - 54 DF_{33}F_3 
+ 186 F_{23}F_3 +\\&& 66 DF_3 F_{33} - 252 F_2 F_{33} - 31 F_3^2
    F_{33})\theta^1_0
 + \frac{1}{90}(12 DF_{33} - 6 F_{23} + F_3 F_{33})\theta^2_0 + \\&&
    \frac{1}{360}(-24 DF_3 + 36 F_2 + 11 F_3^2)\theta^4_0\\
\Om^0&=&\frac{1}{4320}(120 DF_{23} + 240 F_{13} - 240 F_{22} + 36 DF_{33}F_3 - 168 F_{23}F_3 -\\&& 36DF_3 F_{33} + 204 F_2 F_{33} + 17 F_3^2 F_{33})\theta^0_0 + 
    \\&&\frac{1}{12}(-DF_{33} + F23)\theta^1_0 - \frac16 F_{33}\theta^2_0 + \frac{1}{12}F_3\theta^4_0.
\end{eqnarray*}

All the eight forms
$(\theta^0,\theta^1,\theta^2,\theta^3,\Om^0_+,\Om^0_0,\Om^0_-,\Om^0)$
satisfy system (\ref{sy1})-(\ref{sy2}), with corresponding
coefficients $(a_0^0,a_1^0,a_2^0,b_0^0,b_1^0,b_2^0,b_3^0,b_4^0)$ given by:
\begin{eqnarray*}
b_4^0&=&-2F_{333}\\
b_3^0&=&\frac{1}{12}(-6 DF_{333} - 5 F_3 F_{333})\\
b_2^0&=&\frac{1}{360}(120DF_{233} + 240 F_{133} - 120 F_{223} + 60 DF_{333} F_3 - 36DF_3 F_{333} + 
      204 F_2 F_{333} + 79 F_3^2 F_{333})\\
b_1^0&=&\frac{1}{1080}(-180DF_{223} - 540 F_{123} + 360 F_{222} - 90 DF_{33}F_{23} + 
      270 F_{23}^2 + 90 DF_3 F_{233} - 540 F_2 F_{233}\\&& - 180 DF_{233}F_3 - 
      270 F_{133} F_3 + 360 F_{223} F_3 - 45 DF_{333} F_3^2 - 45
      F_{233} 
F_3^2 + 
      90 DF_{23} F_{33} + \\&&90 F_{13}F_{33} - 360 F_{22} F_{33} - 90
      F_{23} 
F_3 F_{33} + 
      90 F_2 F_{33}^2 + 18 DF_2 F_{333} - 72 F_1 F_{333} + \\&&54 DF_3 F_3 F_{333} - 
      288 F_2 F_3 F_{333} - 71 F_3^3F_{333} - 180 F_{33y})
\end{eqnarray*}
\begin{eqnarray*}
b_0^0&=&\frac{1}{129600}(-8640 DF_{233} DF_3 - 12960 DF_{23} DF_{33} - 4320 DF_2 DF_{333} + 
      43200 DF_{33y} + \\&&17280 DF_{333} F_1 + 129600 F_{113} - 64800 F_{122} - 
      34560 DF_{33} F_{13} - 86400 DF_3 F_{133} +\\&& 12960 DF_{233} F_2 + 
      194400 F_{133} F_2 + 30240 DF_{33} F_{22} + 32400 DF_3 F_{223} - 
      64800 F_2 F_{223} +\\&& 6480 DF_{23} F_{23} - 25920 F_{13} F_{23} - 15120 F_{22} F_{23} + 
      2160 DF_2 F_{233} - 8640 F_1 F_{233} - \\&&64800 F_{23y} - 6480 DF_{33}^2 F_3 - 
      6480 DF_3DF_{333} F_3 - 21600 F_{123} F_3 + 10800 DF_{333} F_2 F_3 -\\&& 
      10800 F_{222} F_3 + 25560 DF_{33} F_{23} F_3 - 18360 F_{23}^2 F_3 + 
      13320 DF_3 F_{233} F_3 - 25920F_2F_{233} F_3 +\\&& 3960 DF_{233} F_3^2 + 
      50400 F_{133} F_3^2 - 28800 F_{223} F_3^2 + 2820 DF_{333} F_3^3 - 
      10980 F_{233} F_3^3 - \\&&18000 DF_{22} F_{33} + 6480 DF_3 DF_{33} F_{33} + 
      86400 F_{12} F_{33} - 28080 DF_{33} F_2 F_{33} - \\&&11880 DF_3 F_{23} F_{33} + 
      10800 F_2 F_{23} F_{33} - 19080 DF_{23} F_3 F_{33} + 18720 F_{13} F_3 F_{33} + 
      \\&&16920 F_{22} F_3 F_{33} - 8100 DF_{33} F_3^2 F_{33} + 7200 F_{23} F_3^2 F_{33} + 
      7560 DF_2 F_{33}^2 - 30240 F_1 F_{33}^2 -\\&& 11520 F_2 F_3 F_{33}^2 - 
      1620 F_3^3 F_{33}^2 + 11664 DF_3^2 F_{333} - 63072 DF_3 F_2 F_{333} + 
      76464 F_2^2 F_{333} -\\&& 2520 DF_2 F_3 F_{333} + 10080 F_1 F_3 F_{333} - 
      17712 DF_3 F_3^2 F_{333} + 42768 F_2 F_3^2F_{333} + 5299 F_3^4 F_{333} - \\&&
      18000 F_3 F_{33y} - 75600 F_{33} F_{3y} + 43200 F_{333} F_y)\\
a_2^0&=&\frac{1}{45}(-18 DF_{333} - 24 F_{233} - 4 F_{33}^2 - 27 F_3
F_{333})\\
a_1^0&=&\frac{1}{540}(-72 DF_{233} - 432 F_{133} + 216 F_{223} - 36 DF_{333} F_3 + 
      96 F_{233} F_3 + 48 DF_{33} F_{33} + 16 F_3 F_{33}^2 +\\&& 108 DF_3 F_{333} - 
      324 F_2 F_{333} - 81 F_3^2F_{333})\\
a_0^0&=&\frac{1}{4050}(-180 DF_{223} + 288 DF_{33}^2 - 4860 F_{123} + 2520 F_{222} - 
      378 DF_{33} F_{23} + 1782 F_{23}^2 + \\&&810 DF_3 F_{233} - 2700 F_2 F_{233} - 
      180 DF_{233} F_3 - 2430 F_{133} F_3 + 2880 F_{223} F_3 - 45
      DF_{333} 
F_3^2 + \\&&
      435 F_{233} F_3^2 + 810 DF_{23} F_{33} + 810 F_{13} F_{33} -
      2520 F_{22} 
F_{33} + 
      408 DF_{33} F_3 F_{33} - 594 F_{23} F_3 F_{33} +\\&& 810 F_2 F_{33}^2 + 
      122 F_3^2 F_{33}^2 - 270 DF_2 F_{333} + 1080 F_1 F_{333} + 
      270 DF_3 F_3 F_{333} -\\&& 1080 F_2 F_3 F_{333} - 135 F_3^3 F_{333} + 
      2700 F_{33y}).
\end{eqnarray*}
One can use these, relatively simple, formulae to generate expressions
for the invariant forms on $P$. This may be achieved by means of a matrix
\be
m=\bma
\al^1_{~1}\al^4_{~4}&0&0&0\\
&&&\\
-\frac{\al^1_{~0}}{3}&\al^1_{~1}&0&0\\
&&&\\
\frac{(\al^1_{~0})^2}{9\al^1_{~1}\al^4_{~4}}&-\frac{2\al^1_{~0}}{3\al^4_{~4}}&\frac{\al^1_{~1}}{\al^4_{~4}}&0\\
&&&\\
\frac{-(\al^1_{~0})^3}{27(\al^1_{~1}\al^4_{~4})^2}&\frac{(\al^1_{~0})^2}{3\al^1_{~1}(\al^4_{~4})^2}&-\frac{\al^1_{~0}}{(\al^4_{~4})^2}&\frac{\al^1_{~1}}{(\al^4_{~4})^2}
\ema.
\ee
Then the expression for the invariant 1-forms 
$(\theta^i)=(\theta^0,\theta^1,\theta^2,\theta^3)$ can be
written as
\be
\theta^i=m^i_{~j}\theta^j_0,\quad\quad
i,j=0,1,2,3.
\label{tr1}
\ee
The residual group $G=\{m~|~\al^1_{~1},\al^4_{~4}\neq
0,\al^1_{~0}\in\bbR\}$ has the Lie algebra
$\mathfrak{g}={\mathfrak h}_2\oplus{\mathfrak h}_1$ isomorphic to the
direct sum of the 2-dimensional 
noncommuting Lie algebra ${\mathfrak h}_2$ and a 1-dimensional Lie
algebra ${\mathfrak h}_1$. Algebra ${\mathfrak h}_2$ is 
related to the parameters $(\al^1_{~0},\al^4_{~4})$
and algebra ${\mathfrak h}_1$ is associated
with $\al^1_{~1}$.

The action of $G$ on $\theta^i_0$, induces its action on
$(\Om^0_+,\Om^0_0,\Om^0_-,\Om^0)$. Indeed, defining 
$$\con=\Omega^0_- E_-+\Omega^0_+E_++\Omega^0_0 E_0+\Omega^0
E,$$
and 
$$\Gamma=\Omega_- E_-+\Omega_+E_++\Omega_0 E_0+\Omega E,$$
where $4\times 4$ matrices $(E_-,E_+,E_0,E)$ are the generators
of the Lie algebra $\gla(2,\bbR)$ given in 
(\ref{rep}), we find that
$$\Gamma=m\con m^{-1}+m\der m^{-1}.$$
This enables us to find the explicit expressions for the invariant forms
$(\Om_+,\Om_0,\Om_-,\Om)$.

The transformation rule for $\Gamma$ resembles transformation rule for
a connection. Since $\Gamma$ is $\gla(2,\bbR)$-valued, it is
reasonable to look for a $\glg(2,\bbR)$ principal fibre bundle associated with
corresponding ODE (\ref{ode4}). 

Due to properties of system (\ref{sy1})-(\ref{sy2}) the desired bundle
is just $P$ of Theorem \ref{prp}. To see this, note that equations
  (\ref{sy1}) ensure that $(\theta^1,\theta^2,\theta^3,\theta^4)$ form
  a closed differential ideal. Thus a 4-dimensional distribution
  $\mathcal V$ on $P$ such that ${\mathcal V}\hook\theta^i=0$,
  $\forall i=0,1,2,3$, is integrable. As a consequence, the manifold 
$P$ is foliated by 4-dimensional integral leaves of this
  distribution. Looking at equations (\ref{sy2}) we see that \emph{on each
  leaf} of $\mathcal V$ the forms $(\Om_+,\Om_0,\Om_-,\Om)$ satisfy the Maurer-Caratn
  equations for the $\glg(2,\bbR)$ group. This means that $P$ is a
  principal $\glg(2,\bbR)$ bundle over the leaf space $M^4=P/{\mathcal
  V}$. This 4-dimensional space may be identified with a solutions
  space of ODE (\ref{ode4}). 
\begin{remark}
For local calculations, 
it may be convenient to pass from coordinates 
$(x,y,y_1,y_2,y_3,\al^1_{~0},\al^1_{~1},\al^4_{~4})$ on $P$ to coordinates
$(c_0,c_1,c_2,c_3,s,\al^1_{~0},\al^1_{~1},\al^4_{~4})$ on $P$, where
$(c_0,c_1,c_2,c_3)$ are the integration constants of ODE (\ref{ode4}), 
and $s$ is a real parameter
such that the total differential vector filed $D=\partial_s$. In such
parametrisation $(s,\al^1_{~0},\al^1_{~1},\al^4_{~4})$ constitute
coordinates on the leaves of $\mathcal V$ and 
$(c_0,c_1,c_2,c_3)$ parametrise the solution space $M^4$.   
\end{remark}

\section{$\glg(2,\bbR)$ geometry on the solution space}\label{sec4}
Using matrices $\Gamma=(\Gamma^i_{~j})$, $i,j=0,1,2,3$, and part 
$\theta^i=(\theta^0,\theta^1,\theta^2,\theta^3)$ of the invariant
coframe we rewrite equations (\ref{sy1}) in a compact form as:
\be
\der\theta^i+\Gamma^i_{~j}\dz\theta^j=0,\label{sys1}\ee
and equation (\ref{sy2}) in a compact form as:
\be
\der\Gamma^i_{~k}+\Gamma^i_{~j}\dz\Gamma^j_{~k}=\tfrac12 R^i_{~kjl}\theta^j\dz\theta^l.\label{sys2}\ee
The coefficients $R^i_{~jkl}$ appearing in this last equation can be
easily read off from (\ref{sy2}). They are linear combinations of the
coefficients $a_0,a_1,a_2,a_3,b_0,b_1,b_2,b_3,b_4$ of
(\ref{sy2}). The meaning of equations (\ref{sys1})-(\ref{sys2}) is
obvious: they constitute, respectively, the first and the second
Cartan's structure equations, for a $\gla(2,\bbR)$-valued
connection $\Gamma$ on the principal fibre bundle $\glg(2,\bbR)\to P\to
M^4$. Due to the first
equation, (\ref{sys1}), this connection has no torsion. The second
equation, (\ref{sys2}), determines the curvature of $\Gamma$; the
coefficients $R^i_{~jkl}$ are the curvature tensor coefficients for
$\Gamma$. 

Given the curvature tensor $R^i_{~jkl}$ of $\Gamma$ we define 
its `Ricci' tensor $R_{jl}$ by
$$R_{jl}=R^i_{~jil}.$$

Recalling that the curvature of $\Gamma$ is totally expressible in terms 
of $a_0,a_1,a_2,a_3,b_0,$ $b_1,b_2,b_3,b_4$ and performing 
a purely algebraic manipulation
on the curvature tensor coefficients $R^i_{~jkl}$ we get a remarkable 
\begin{theorem}\label{there}
Every 4th order ODE satisfying conditions (\ref{wun}) uniquely defines
a principal fibre bundle $\glg(2,\bbR)\to P\to M^4$ over the space of
its solutions $M^4$ and a torsionless $\gla(2,\bbR)$-connection
$\Gamma$ on $P$ with curvature $R^i_{~jkl}$ having the Ricci tensor
$R_{jl}$ in the form
$$
R_{jl}=\bma 0&b_0&a_0+2b_1&-a_1+b_2\\
-b_0&-2a_0&a_1+3b_2&a_2+2b_3\\
a_0-2b_1&a_1-3b_2&-2a_2&b_4\\
-a_1-b_2&a_2-2b_3&-b_4&0
\ema.
$$
Its respective symmetric and antisymmetric parts read:
$$
R_{(jl)}=\bma 0&0&a_0&-a_1\\
0&-2a_0&a_1&a_2\\
a_0&a_1&-2a_2&0\\
-a_1&a_2&0&0
\ema,
$$
and 
$$
R_{[jl]}=\bma 0&b_0&2b_1&b_2\\
-b_0&0&3b_2&2b_3\\
-2b_1&-3b_2&0&b_4\\
-b_2&-2b_3&-b_4&0
\ema.
$$
Thus the entire curvature tensor $R^i_{~jkl}$ is encoded in the Ricci tensor.
\end{theorem}
\begin{remark}
Note that we also have $R^i_{~ikl}=2R_{[kl]}$.
\end{remark}

Now we can use matrix $m$ of the previous section to find explicit
formulae for the coefficients
$a_0,a_1,a_2,a_3,b_0,b_1,b_2,b_3,b_4$. It follows that if we evaluate 
$R_{ij}$ for $(\al^1_{~0},\al^1_{~1},\al^4_{4})=(0,1,1)$, 
denoting the calculated $R_{ij}$ by $R^0_{ij}$, then the full Ricci
tensor $R_{ij}$ is related to $R^0_{ij}$ via
\be
R_{ij}=R^0_{kl}{m^{-1}}^k_{~i}{m^{-1}}^l_{~j}. 
\label{tr2}
\ee
Here $m^{-1}=({m^{-1}}^j_{~j})$ is the inverse matrix to $m$. From this expression we
can calculate explicit form of
$a_0,a_1,a_2,a_3,b_0,b_1,b_2,b_3,b_4$. The resulting formulae involve 
coefficients $a^0_0,a^0_1,a^0_2,a^0_3,b^0_0,b^0_1,b^0_2,b^0_3,b^0_4$ of
the previous section and parameters
$\al^1_{~0},\al^1_{~1},\al^4_{4}$ and read:
\begin{eqnarray*}
b_4&=&\frac{(\al^4_{~4})^3}{(\al^1_{~1})^2}b^0_4\\
b_3&=&\frac{(\al^4_{~4})^2}{(\al^1_{~1})^2}b^0_3 +
\frac{\al^1_{~0}(\al^4_{~4})^2}{3(\al^1_{~1})^3}b^0_4\\
b_2&=&\frac{\al^4_{~4}}{(\al^1_{~1})^2}b^0_2 + \frac{2\al^1_{~0}
  \al^4_{~4}}{3(\al^1_{~1})^3}b^0_3 +
\frac{(\al^1_{~0})^2\al^4_{~4}}{9(\al^1_{~1})^4}b^0_4\\
b_1&=&\frac{1}{(\al^1_{~1})^2}b^0_1 +
\frac{\al^1_{~0}}{(\al^1_{~1})^3}b^0_2
+\frac{(\al^1_{~0})^2}{(3\al^1_{~1})^4}b^0_3
+\frac{(\al^1_{~0})^3}{27(\al^1_{~1})^5}b^0_4\\
b_0&=&\frac{1}{(\al^1_{~1})^2\al^4_{~4}}b^0_0 +\frac{4\al^1_{~0}}{3(\al^1_{~1})^3\al^4_{~4}}b^0_1 +\frac{2(\al^1_{~0})^2}{3(\al^1_{~1})^4\al^4_{~4}}b^0_2 + \frac{4(\al^1_{~0})^3}{27(\al^1_{~1})^5\al^4_{~4}}b^0_3 + \frac{(\al^1_{~0})^4}{81(\al^1_{~1})^6\al^4_{~4}}b^0_4,
\end{eqnarray*}
\begin{eqnarray}
a_2&=&\frac{(\al^4_{~4})^2}{(\al^1_{~1})^2}a^0_2\nonumber\\
a_1&=&\frac{\al^4_{~4}}{(\al^1_{~1})^2}a^0_1 -
\frac{\al^1_{~0}\al^4_{~4}}
{3(\al^1_{~1})^3}a^0_2\label{atr}\\
a_0&=&\frac{1}{(\al^1_{~1})^2}a^0_0 -\frac{2\al^1_{~0}}{3(\al^1_{~1})^3}a^0_1 +\frac{(\al^1_{~0})^2}{9(\al^1_{~1})^4}a^0_2.\nonumber
\end{eqnarray}
This, in particular, means that the respective spaces
consisting of $(b^0_0,b^0_1,b^0_2,b^0_3,b^0_4)$ and of 
$(a^0_0,a^0_1,a^0_2)$ constitute
a 5-dimensional and 3-dimensional representation of $G$ and, as a
consequence of $\glg(2,\bbR)$.

Due to (\ref{tr2}) the vanishing of any of the two determinants: 
$$\det( R_{(ij)})\quad\quad{\rm and}\quad\quad \det(
R_{[ij]})$$ 
is a contact invariant property of the
corresponding 4th order ODE (\ref{ode}). These two determinants, when
expressed in terms of the eight curvature coefficients
$a_0,a_1,a_2,b_0$, $b_1,b_2,b_3,b_4$, are 
$$\det( R_{(ij)})=(a_1^2-a_0a_2)^2$$
and
$$\det( R_{[ij]})=(3b_2^2-4b_1b_3+b_0b_4)^2.$$
Thus they are expressible in terms of the two well known 
$\glg(2,\bbR)$-invariant polynomials
$$I_2=a_1^2-a_0a_2\quad\quad{\rm and}\quad\quad
I_3=3b_2^2-4b_1b_3+b_0b_4.$$
\begin{remark}
In this context it is interesting to note that function
$F=(y_3)^{(4/3)}$ of the well known example (\ref{examm}), provides 
a contact equivalent class of ODEs that has both 
invariants $I_2$ and $I_3$ vanishing.
\end{remark}

Interestingly, the next $\glg(2,\bbR)$-invariant polynomial
$$I_4=-3(\theta^1)^2(\theta^2)^2 + 4\theta^0(\theta^2)^3 + 4(\theta^1)^3 \theta^3 - 6\theta^0\theta^1\theta^2\theta^3 + (\theta^0)^2(\theta^3)^2,$$
when thought as defined on $P$ in terms of forms
$(\theta^0,\theta^1,\theta^2,\theta^3)$ of the invariant
coframe $(\theta^0,\theta^1,\theta^2,\theta^3,\Om_+,\Om_-,\Om_0,\Om)$,
has the following property:
$${\mathcal L}_X I_4=12(X\hook\Om) I_4,$$
where $X\in{\mathcal V}$ is any vertical vector field on $\glg(2,\bbR)\to P\to
M^4$. Thus
$I_4$ descends to a well defined conformal
symmetric tensor of fourth degree on the solution space $M^4$ of the
ODE \cite{bryantsp4}. Let us denote the descended to $M^4$ tensor
$I_4$ by $\ten$. It is also worthwhile to mention that, for the vertical vectors
$X\in\mathcal V$, we have  
$${\mathcal L}_X\Om=\der(X\hook\Om).$$
This means that on the solution space $M^4$ the form $\Om$ is defined up to
a \emph{gradient}. It is convenient to rescale $\Om$ and to define a 1-form
$A$ on $P$ equal to $$A=-12\Om.$$ This form is also defined up to a
gradient on the solutions space $M^4$. 
Thus, a solution space $M^4$ of any 4th order ODE satisfying
(\ref{wun}) is equipped with a sort of Weyl geometry $[\ten,A]$. This
consists of class of pairs $(\ten,A)$, in which $\ten$ is a 4th order
symmetric tensor field, $A$ is a 1-form on $M^4$, and  two pairs
$(\ten,A)$ and $(\ten',A')$ represent the same class iff
$$\ten'={\rm e}^{4\phi}\ten,\quad\quad A'=A-4\der\phi.$$
In the context of this gauge freedom, 
it is worthwhile to note that the vanishing of
$R_{[ij]}$ corresponds to the $[\ten,A]$ geometries on $M^4$ with form
$A$ that can be gauged to $A=0$. Such situation occurs if and only if 
$b_i=0$ for all $i=0,1,2,3,4$.

\begin{remark}
In terms of the Weyl-like 
geometry $[\ten,A]$ on the solution space $M^4$, 
the $\gla(2,\bbR)$-valued connection may be
defined as the unique torsionless connection satisfying
$$\nabla_X \ten=-A(X)\ten.$$
\end{remark}          

Thus we have the following
\begin{theorem}
Every 4th order ODE $y^{(4)}=F(x,y,y',y'',y^{(3)})$ satisfying
Braynt's conditions (\ref{wun}) uniquely defines a conformal Weyl-like
geometry
$[\ten,A]$ on its solution space $M^4$. The Weyl-like geometry
$[\ten,A]$ 
consists of a symmetric 4th rank tensor $\ten$ and a
1-form $A$ given up to transformations  
$$\ten'={\rm e}^{4\phi}\ten,\quad\quad A'=A-4\der\phi.$$ Its
corresponding $\gla(2,\bbR)$-valued connection has no torsion and very
special curvature tensor described by Theorem \ref{there}.
\end{theorem}
\section{Examples}
\subsection{Equations with symmetric Ricci tensor}
There is only one contact equivalence class of ODEs (\ref{ode4})
having an 8-dimensional group of contact symmetries. This is
equivalent to $y^{(4)}=0$ and the symmetry group is
$\glg(2,\bbR)\times_{\rho}\bbR^4$. For this class of equation the
$\gla(2,\bbR)$-valued connection of Theorem \ref{there} is flat.

In this section we focus on the equivalence classes of ODEs
(\ref{ode4}) for which the Maxwell form $\der A=-12\der\Om$ of this 
connection is flat $\der F=0$. In such case we have
$b_0=b_1=b_2=b_3=b_4=0$.

Let us assume that we are in this situation.

Looking at the transformation properties (\ref{atr}) of the curvature
coefficient $a_2$ we see that there are essentially two distinct
cases distinguished by the \emph{vanishing or not} of the expression 
$a^0_2=\frac{1}{45}(-18 DF_{333} - 24 F_{233} - 4 F_{33}^2 - 27
F_3F_{333})$. 

We analyse the more easy case $a^0_2=0$ first.

If $a^0_2=0$ then also $a_2=0$. Thus we have $a_2=0$ everywhere on $P$
  with the full system (\ref{sy1})-(\ref{sy2}) of eight independent
  1-forms $\theta^1,\theta^2,\theta^3,\Om_+,\Om_0,\Om_-,\Om$
  there. Imposing $(\der^2\Om_+)\dz\theta^1\dz\theta^2=0$ on 
(\ref{sy1})-(\ref{sy2}) quickly leads to $a_1=0$ and, consequently, by imposition of
  $(\der^2\Om_+)\dz\theta^1=0$, to
  $a_0=0$. This shows that if $a^0_2=0$ then the
  corresponding ODEs (\ref{ode4}) are contact equivalent to
  $y^{(4)}=0$.

Now we assume that $a^0_2\neq 0$. Then the choice $$\al^1_{~1}=\frac{\sqrt{2}}{4}\al^4_{~4}\sqrt{|a^0_2|}$$
brings $a_2$ to the form $$a_2=8\epsilon_1,$$ where
$\epsilon_1=\sgn(a^0_2)$. Then the choice
$$\al^1_{~0}=\frac{3\sqrt{2}}{4}\epsilon_1\al^4_{4}\frac{a^0_1}{\sqrt{|a^0_2|}}$$ makes
$$a_1=0.$$ After these two normalisations we get 
$$a_0=\frac{8\epsilon_1}{(\al^4_{~4}a^0_2)^2)}(a^0_0a^0_2-(a^0_1)^2).$$
Thus again we have two cases, depending on the vanishing or not of
the invariant $I^0_2=(a^0_1)^2-a^0_0a^0_2$.

It follows that the $I^0_2=0$ case, which under our assumptions is the same as $a_0=0$,
corresponds to only \emph{one} nonequivalent class of equations. 
They are defined by $\epsilon_1=1$ (the
$\epsilon_1=-1$ case is not compatible with system (\ref{sy1})-(\ref{sy2})), and are described by the following
\begin{theorem}\label{street}
All ODEs $y^{(4)}=F(x,y,y',y'',y^{(3)})$ satisfying Bryant's conditions
(\ref{wun}), having symmetric Ricci tensor, and invariants 
$I_2=0$ and $a_2\neq 0$, are in local one-to-one 
correspondence with coframes 
$(\theta^0,\theta^1,\theta^2,\theta^3,\Om_+,\Om)$ on a 6-manifold
satisfying:
\begin{eqnarray*}
\der\theta^0&=&12\Om\dz\theta^0-3\Om_+\dz\theta^1+\tfrac{3\sqrt{2}}{2}\theta^0\dz\theta^2\\
\der\theta^1&=&6\Om\dz\theta^1-2\Om_+\dz\theta^2+\tfrac{\sqrt{2}}{2}(\theta^0\dz\theta^3+\theta^1\dz\theta^2)\\
\der\theta^2&=&-\Om_+\dz\theta^3+\sqrt{2}\theta^1\dz\theta^3\\
\der\theta^3&=&-6\Om\dz\theta^3+ 3\sqrt{2}\theta^2\dz\theta^3\\
\der\Om_+&=&6\Om\dz\Om_++\sqrt{2}\Om_+\dz\theta^2+\theta^0\dz\theta^3-5\theta^1\dz\theta^2\\
\der\Om&=&0.
\end{eqnarray*}
The forms $\Om$ and $\Om_0$ are given by
$$\Om_-=\frac{\sqrt{2}}{2}\theta^3,\quad\quad
\Om_0=3\Om-\frac{\sqrt{2}}{2}\theta^2.$$
All the equations having such invariant forms are equivalent to
an ODE defined by 
$$F=\frac43 \frac{y_3^2}{y_2}.$$
This class has strictly 6-dimensional group of contact symmetries.
\end{theorem}

Now we pass to the $I_2^0\neq 0$ case. We introduce $\epsilon_2=\pm 1$, which encodes the sign of
$I^0_2$. This is defined by  
$\epsilon_1\epsilon_2(a^0_0a^0_2-(a^0_1)^2)> 0$. Now we 
chose
$$\al^4_{~4}=\sqrt{\frac{\epsilon_1\epsilon_2(a^0_0a^0_2-(a^0_1)^2)}{(a^0_2)^2}}.$$
This normalises $a_0$ to
$$a_0=8\epsilon_2.$$

Under such normalisations system (\ref{sy1})-(\ref{sy2})
descends from $P$ to the 5-dimensional jet space $J$. There, it reads:
\begin{eqnarray}
\der\theta^0&=&3(\Om+\Om_0)\dz\theta^0-3\Om_+\dz\theta^1\nonumber\\
\der\theta^1&=&-\Om_-\dz\theta^0+(3\Om+\Om_0)\dz\theta^1-2\Om_+\dz\theta^2\nonumber\\
\der\theta^2&=&-2\Om_-\dz\theta^1+(3\Om-\Om_0)\dz\theta^2-\Om_+\dz\theta^3\nonumber\\
\der\theta^3&=&-3\Om_-\dz\theta^2+3(\Om-\Om_0)\dz\theta^3\nonumber\\&&\nonumber\\
\der\Om_+&=&2\Om_0\dz\Om_+-
2\epsilon_2 \theta^0\dz\theta^1 +  
    \epsilon_1(\theta^0\dz\theta^3 -
    5 \theta^1\dz\theta^2) \nonumber\\
\der\Om_-&=&-2\Om_0\dz\Om_-+
    \epsilon_2 (-\theta^0\dz\theta^3 + 
    5\theta^1\dz\theta^2)+  
    2\epsilon_1\theta^2\dz\theta^3\nonumber\\
\der\Om_0&=&\Om_+\dz\Om_- -
    \epsilon_2 \theta^0\dz\theta^2 - 
    \epsilon_1 \theta^1\dz\theta^3 \nonumber\\
\der\Om&=&0\nonumber
\end{eqnarray} 
To close this system it is convenient to eliminate form $\Omega$. This
can be achieved by an introduction of new forms
$(\sigma^0,\sigma^1,\sigma^2,\sigma^3)$ related to
$(\theta^0,\theta^1,\theta^2,\theta^3)$ via:
$$\sigma^0={\rm e}^{w}\theta^0,\quad\sigma^1={\rm
  e}^{w}\theta^1,\quad\sigma^2={\rm
  e}^{w}\theta^2,\quad\sigma^3={\rm e}^{w}\theta^3,$$
where $w$ is a function on $J$ such that $\Om=-\tfrac13\der w$. The
  local existence of such function is guaranteed by $\der\Om=0$. In
  terms of the new variables $(\sigma^0,\sigma^1,\sigma^2,\sigma^3)$,
  $w$, the reduced system takes a form in which the 1-form $\Om$ is
  not present:
\begin{eqnarray}
\der\sigma^0&=&3\Om_0\dz\sigma^0-3\Om_+\dz\sigma^1\nonumber\\
\der\sigma^1&=&-\Om_-\dz\sigma^0+\Om_0\dz\sigma^1-2\Om_+\dz\sigma^2\nonumber\\
\der\sigma^2&=&-2\Om_-\dz\sigma^1-\Om_0\dz\sigma^2-\Om_+\dz\sigma^3\nonumber\\
\der\sigma^3&=&-3\Om_-\dz\sigma^2-3\Om_0\dz\sigma^3\nonumber\\&&\label{sya}\\
\der\Om_+&=&2\Om_0\dz\Om_+ +{\rm e}^{-2w}\Big(-
2\epsilon_2\sigma^0\dz\sigma^1 +  
    \epsilon_1(\sigma^0\dz\sigma^3 -
    5 \sigma^1\dz\sigma^2)\Big) \nonumber\\
\der\Om_-&=&-2\Om_0\dz\Om_-+{\rm e}^{-2w}\Big(
    \epsilon_2 (-\sigma^0\dz\sigma^3 + 
    5\sigma^1\dz\sigma^2)+  
    2\epsilon_1\sigma^2\dz\sigma^3\Big)\nonumber\\
\der\Om_0&=&\Om_+\dz\Om_- -{\rm e}^{-2w}\Big(
    \epsilon_2 \sigma^0\dz\sigma^2 + 
    \epsilon_1 \sigma^1\dz\sigma^3\Big). \nonumber
\end{eqnarray} 
As we can see the price paid for elimination of $\Om$ is and
introduction of nonconstant function $w$ appearing explicitly in
these equations.

Now the remarkable fact is that system (\ref{sya})
closes on $J$ and is described by the following Theorem.
\begin{theorem}
All ODEs $y^{(4)}=F(x,y,y',y'',y^{(3)})$ satisfying Bryant's conditions
(\ref{wun}), having symmetric Ricci tensor, and invariants 
$I_2\neq 0$ and $a_2\neq 0$, are in local one-to-one 
correspondence with coframes 
$(\sigma^0,\sigma^1,\sigma^2,\sigma^3,\Om_+)$ on a 5-manifold
satisfying system (\ref{sya}) with:
\begin{eqnarray}
\Om_0&=&w_0\sigma^0-(w_1+4\epsilon_1\epsilon_2w_3)\sigma^1+(4\epsilon_1\epsilon_2w_0+w_2)\sigma^2-w_3\sigma^3\nonumber\\
\Om_-&=&-\epsilon_1\epsilon_2\Om_+-2(\epsilon_1\epsilon_2w_1+2w_3)\sigma^0+2w_0\sigma_1+\label{syb}\\&&2\epsilon_1\epsilon_2w_3\sigma^2-2(2\epsilon_1\epsilon_2w_0+w_2)\sigma^3.\nonumber
\end{eqnarray}
Functions $w,w_0,w_1,w_2,w_3$ appearing here are defined by:
\be
\der w=w_0\sigma^0+w_1\sigma^1+w_2\sigma^2+w_3\sigma^3.
\label{syc}\ee
They satisfy
\begin{eqnarray}
\der w_0 &=& -\epsilon_1 \epsilon_2 w_1 \Om_+ +\frac14(-\epsilon_2{\rm
   e}^{-2w}+4w_0^2+16 w_1 w_3+32\epsilon_1\epsilon_2 w_3^2)\sigma^0 +\nonumber\\&&
   3 w_0 w_1 \sigma^1
     - (-\epsilon_1 \epsilon_2 w_{13} - 11 w_0 w_2 - 4 \epsilon_1 \epsilon_2 w_2^2 +
      5 \epsilon_1 \epsilon_2 w_1 w_3 + 12 w_3^2) \sigma^2 +\nonumber\\&&(11 w_0 + 4 \epsilon_1 \epsilon_2 w_2) w_3 \sigma^3\nonumber\\
\der w_1& =& (3 w_0 - 2 \epsilon_1 \epsilon_2 w_2) \Om_+ + 
(-3 w_0 w_1 - 4 \epsilon_1 \epsilon_2 w_1 w_2 -
     12 \epsilon_1 \epsilon_2 w_0 w_3 - 8 w_2 w_3) \sigma^0 -\nonumber\\&&
   \frac14(3 \epsilon_1{\rm e}^{-2w} +
     24 \epsilon_1 \epsilon_2  w_0^2 -
     20 w_1^2 +
     8 \epsilon_1 \epsilon_2 w_{13} +
     64  w_0 w_2 +
     32 \epsilon_1 \epsilon_2 w_2^2 -\nonumber\\&&
     120 \epsilon_1 \epsilon_2 w_1 w_3 -
     192 w_3^2) \sigma^1 - \nonumber\\&&(12 w_0 w_1 \epsilon_1 \epsilon_2 + w_1 w_2
   + 30 w_0 w_3 + 4 \epsilon_1 \epsilon_2 w_2 w_3) \sigma^2 + w_{13}
   \sigma^3\label{syd}
\end{eqnarray}
\begin{eqnarray*}
\der w_2& = &(2 w_1 -
     3 \epsilon_1 \epsilon_2 w_3) \Om_+ +\\&& \frac12(24 \epsilon_1 \epsilon_2 w_0^2 + 2 \epsilon_1 \epsilon_2 w_{13} +
     30  w_0 w_2 +
     8  \epsilon_1 \epsilon_2 w_2^2 -
     26  \epsilon_1 \epsilon_2 w_1 w_3 -
     48  w_3^2) \sigma^0 - \\&&(8 \epsilon_1 \epsilon_2 w_0 w_1 + w_1 w_2 + 24  w_0 w_3 + 12 \epsilon_1 \epsilon_2 w_2 w_3) \sigma^1 +\\&&
  \frac14 (-3\epsilon_2 {\rm e}^{-2w} +
     96  w_0^2 - 8  w_{13} -
     12  w_2^2 + 40  w_1 w_3 +
     96 \epsilon_1 \epsilon_2 w_3^2) \sigma^2 -\\&&
  3 (8 \epsilon_1 \epsilon_2 w_0 + 3 w_2) w_3 \sigma^3\\
\der w_3 &=& w_2 \Om_+ + (4 \epsilon_1 \epsilon_2 w_0 w_1 + 2 w_1 w_2 + 11 w_0 w_3 +
     4 \epsilon_1 \epsilon_2 w_2 w_3) \sigma^0 + \\&&(w_{13} + 8 \epsilon_1 \epsilon_2 w_0 w_2 + 4 w_2^2 - 4 w_1 w_3 -
     12 \epsilon_1 \epsilon_2 w_3^2) \sigma^1 +
  w_2 w_3 \sigma^2 +\\&& \frac14 (-\epsilon_1{\rm e}^{-2w} +
     32  \epsilon_1 \epsilon_2 w_0^2 +
     32 w_0 w_2 +
     8 \epsilon_1 \epsilon_2 w_2^2 +
     4 w_3^2) \sigma^3,
\end{eqnarray*}
with a function $w_{13}$ satisfying
\begin{eqnarray}
\der w_{13}&=& (-12 \epsilon_1 \epsilon_2 w_0 w_1 - w_1 w_2 + 45 w_0 w_3 + 
     30 \epsilon_1 \epsilon_2 w_2 w_3) \Omega_+ +\nonumber\\&& 
  \frac12(-6 \epsilon_2 w_0{\rm e}^{-2w} - 
     240  w_0^3 + 
     40  \epsilon_1 \epsilon_2 w_0 w_1^2 - 
     16  w_0 w_{13} + 5 \epsilon_1 w_2{\rm e}^{-2w} - \nonumber\\&&
     472  \epsilon_1 \epsilon_2 w_0^2 w_2 + 
     20 w_1^2 w_2 - 
     16  \epsilon_1 \epsilon_2 w_{13} w_2 - 
     304  w_0 w_2^2 - 
     64  \epsilon_1 \epsilon_2 w_2^3 + \nonumber\\&&
     384  w_0 w_1 w_3 + 
     192  \epsilon_1 \epsilon_2 w_1 w_2 w_3 + 
     552 \epsilon_1 \epsilon_2 w_0 w_3^2 + 
     272  w_2 w_3^2) \sigma^0 - \nonumber\\&&
  \frac14 (20 \epsilon_2 w_1{\rm e}^{-2w} - 
       256  w_0^2 w_1 - 
       28  w_1 w_{13} - 
       416  \epsilon_1 \epsilon_2 w_0 w_1 w_2 - 
       144 w_1 w_2^2 +\nonumber\\&& 15 \epsilon_1  w_3 {\rm e}^{-2w}- 
       840  \epsilon_1 \epsilon_2 w_0^2 w_3 + 
       20  w_1^2 w_3 - 
       24 \epsilon_1 \epsilon_2 w_{13} w_3 - 
       1440 w_0 w_2 w_3 -\nonumber\\&& 
       480 \epsilon_1 \epsilon_2 w_2^2 w_3 - 
       40 \epsilon_1 \epsilon_2 w_1 w_3^2 - 
       192  w_3^3) \sigma^1 -\label{sye}\\&& 
  \frac12(-15 \epsilon_1 w_0{\rm e}^{-2w} + 
       480 \epsilon_1 \epsilon_2 w_0^3 - 
       24 \epsilon_1 \epsilon_2 w_0 w_{13} - 2 \epsilon_2 w_2{\rm e}^{-2w} + 
       544 w_0^2 w_2 - \nonumber\\&&
       16 w_{13} w_2 + 
       184 \epsilon_1 \epsilon_2 w_0 w_2^2 + 
       16 w_2^3 + 
       240  \epsilon_1 \epsilon_2 w_0 w_1 w_3 + 
       80  w_1 w_2 w_3 +\nonumber\\&& 
       588  w_0 w_3^2 + 
       184  \epsilon_1 \epsilon_2 w_2 w_3^2) \sigma^2 -\nonumber\\&& 
  \frac14(5 \epsilon_1 w_1{\rm e}^{-2w} - 
     160 \epsilon_1 \epsilon_2 w_0^2 w_1 - 
     160w_0 w_1 w_2 - 
     40  \epsilon_1 \epsilon_2 w_1 w_2^2 + \nonumber\\&&36 \epsilon_2 w_3{\rm e}^{-2w} - 
     1152  w_0^2 w_3 - 
     28  w_{13} w_3 - 
     1152  \epsilon_1 \epsilon_2 w_0 w_2 w_3 -\nonumber\\&& 
     288  w_2^2 w_3 + 
     20 w_1 w_3^2) \sigma^3.\nonumber
\end{eqnarray}
System (\ref{sya})-(\ref{sye}) is closed, meaning that $\der^2=0$ does
not implies any further relations between forms
$\sigma^0,\sigma^1,\sigma^2,\sigma^3,\Om_+$ and functions $w,w_0,w_1,w_2,w_3,w_{13}$. 
\end{theorem}
We easily see that the assumption that \emph{all}
$w,w_0,w_1,w_2,w_3,w_{13}$ are \emph{constant} is incompatible with
system (\ref{sya})-(\ref{sye}). Finding \emph{any} solution to system
(\ref{sya})-(\ref{sye}) is a difficult tusk. 
\subsection{Inhomogeneous examples}
Here we present examples of contact equivalent classes
of 4th order ODEs satisfying Bryant's conditions (\ref{wun}) which
are \emph{not homogeneous}. By this we mean they do \emph{not} admit a transitive
contact symmetry group of dimension greater than \emph{four}. 
We consider an ansatz in which function
$F$ depends in a special way on only two coordinates $y_2$ and $y_3$. 
Explicitly:
\be
F=~(y_2)^{2}~q\Big(\frac{y_3^{2}}{y_2^{3}}\Big),\label{gor}
\ee
where $q=q(z)$ is a sufficiently differentiable real function of its
argument $$z=\frac{y_3^2}{y_2^3}.$$ 
Imposing Bryant's conditions (\ref{wun}) on (\ref{gor}) we find the
following
\begin{proposition}
Function $F$ of (\ref{gor}) satisfies Bryant's conditions (\ref{wun})
if and only if 
\begin{itemize}
\item[a)] either: $$6z(3z-2q)q''+3z{q'}^2-6qq'+4q=0,$$
\item[b)] or: $$6z(3z-2q)q''+3z{q'}^2-6qq'+14q-15z=0.$$
\end{itemize}
\end{proposition}
The special solutions of a) are: $q(z)=0$ and $q(z)=\frac43 z$. In case
b) we have $q(z)=3z$ and $q(z)=\frac53 z$ as special
solutions. Writing these four
solutions as $q(z)=c z$ we remark that in cases $c=0$ and $c=3$  
function $F$ defines a 4th order ODE which is contact equivalent to
$y^{(4)}=0$. Cases $c=\frac43$ and $c=\frac53$ define two different
$F$s, but the corresponding 4th order ODEs are contact
equivalent. They both are equivalent to the 
ODE described by Theorem \ref{street}. 

We emphasise that apart from the singular solutions $q=cz$, each
equation a) or b) admits a 2-parameter family of solutions. Every
solution $q=q(z)$ from these two families leads to a 4th order ODE which
satisfies Bryant's conditions (\ref{wun}) and which is
\emph{inhomogeneous}. Remarkably all Bryant's $F$s which are defined by the
ansatz (\ref{gor}) have $I_3=I_4=0$, but
$a_2\neq 0$ and $b_4\neq 0$. Thus, in particular, $\der A\neq 0$ for them. 

We were unable to find any example of Bryant's ODEs for which at least
one of $I_2$ or $I_3$ is not vanishing.

\end{document}